# Estimation of volatility functionals in the simultaneous presence of microstructure noise and jumps

MARK PODOLSKIJ[1] and MATHIAS VETTER[2]

[1] *University of Aarhus and CREATES, Department of Economics and Management, 8000 Aarhus C, Denmark.* E-mail: mpodolskij@creates.au.dk
[2] *Ruhr-Universität Bochum, Fakultät für Mathematik, 44780 Bochum, Germany.* E-mail: mathias.vetter@rub.de

We propose a new concept of modulated bipower variation for diffusion models with microstructure noise. We show that this method provides simple estimates for such important quantities as integrated volatility or integrated quarticity. Under mild conditions the consistency of modulated bipower variation is proven. Under further assumptions we prove stable convergence of our estimates with the optimal rate $n^{-1/4}$. Moreover, we construct estimates which are robust to finite activity jumps.

## 1. Introduction

Continuous time stochastic models represent a widely accepted class of processes in mathematical finance. Itô diffusions, which are characterised by the equation

$$X_t = X_0 + \int_0^t a_s \, \mathrm{d}s + \int_0^t \sigma_s \, \mathrm{d}W_s, \qquad (1.1)$$

are commonly used for modelling the dynamics of interest rates or stock prices. Here $W$ denotes a Brownian motion, $a$ is a locally bounded predictable drift function and $\sigma$ is a cadlag volatility process. A key issue in econometrics is the estimation (and forecasting) of the quadratic variation of $X$

$$IV = \int_0^1 \sigma_s^2 \, \mathrm{d}s,$$

which is known as integrated volatility or integrated variance in the econometric literature. In recent years the availability of high frequency data on financial markets has motivated a huge number of publications devoted to measurement of the integrated volatility. A typical way to estimate the integrated volatility is to use the realised volatility (RV), which has been proposed by Andersen, Bollerslev, Diebold and Labys [3] and Barndorff-Nielsen and Shephard [7]. RV is the sum of squared increments over non-overlapping intervals within a sampling period. The consistency result justifying this estimator is a simple consequence of the definition of the quadratic variation (see, e.g., Protter [21]). Theoretical and empirical properties of the realised volatility have been studied in numerous articles (see Jacod [17]; Jacod and Protter [19]; Andersen, Bollerslev, Diebold and Labys [3]; Barndorff-Nielsen and Shephard [7] among many others).

More recently, the concept of realised bipower variation has built a nonparametric framework for backing out several variational measures of volatility (see, e.g., Barndorff-Nielsen and Shephard [8] or Barndorff-Nielsen, Graversen, Jacod, Podolskij and Shephard [5]), which has led to a new development in econometrics. Realised bipower variation, which is defined by

$$BV(X, r, l)_n = n^{(r+l)/2-1} \sum_{i=1}^{n-1} |\Delta_i^n X|^r |\Delta_{i+1}^n X|^l, \qquad (1.2)$$

with $\Delta_i^n X = X_{i/n} - X_{(i-1)/n}$ and $r, l \geq 0$, provides a whole class of estimators for different (integrated) powers of volatility. Another important feature of realised bipower variation is its robustness to finite activity jumps when estimating the integrated volatility (in the case $r \vee l < 2$). This property has been used to construct tests for jumps (see Barndorff-Nielsen and Shephard [9] or Christensen and Podolskij [11]).

However, in finance it is widely accepted that the true price process is contaminated by microstructure effects, such as price discreteness or bid-ask spreads, among others. This invalidates the asymptotic properties of RV, and in the presence of microstructure noise RV is both biased and inconsistent (see Bandi and Russell [4] or Hansen and Lunde [15] among others). Nowadays there exist two concurrent methods of estimating the integrated volatility in the presence of i.i.d. noise. Zhang [24] has proposed to use a *multiscale* estimator as a generalisation of the concept of *two scale* estimators, which was introduced by Zhang, Mykland and Ait-Sahalia [25] based on a subsampling procedure. Another method is a *realised kernel* estimator, which has been proposed by Barndorff-Nielsen, Hansen, Lunde and Shephard [6]. Both methods provide consistent estimates of the integrated volatility in the presence of i.i.d. noise and achieve the optimal rate $n^{-1/4}$ (whereas the two scale approach achieves the rate $n^{-1/6}$). However, these procedures can not be generalised in an obvious way in order to obtain estimators of other (integrated) powers of volatility, such as the integrated quarticity, which is defined by

$$IQ = \int_0^1 \sigma_s^4 \, ds.$$

This quantity is of particular interest because, properly scaled, it occurs as the conditional variance in central limit theorems for estimators of $IV$ and has therefore to be estimated.

Moreover, both methods are not robust to jumps in the price process (here we would like to mention the work by Fan and Wang [12], who obtain jump-robust estimates of *IV* by applying wavelet methods).

In this paper we propose a new concept of modulated bipower variation (MBV) for diffusion models with (i.i.d.) microstructure noise. The novelty of this concept is twofold. First, this method provides a whole class of estimates for arbitrary integrated powers of volatility. Second, modulated multipower variation, which is a direct generalisation of MBV, turns out to be robust to finite activity jumps (when the powers are appropriately chosen). In particular, starting with MBV we construct estimators of *IV* and *IQ* which are robust to finite activity jumps. An easy implementation of MBV is another nice feature of our method.

This paper is organised as follows: In Section 2 we state the basic notation and definitions. In Section 3 we show the consistency of our estimators and prove a central limit theorem for their normalised versions with an optimal rate $n^{-1/4}$. In particular, we construct some new estimators of the integrated volatility and the integrated quarticity and present the corresponding asymptotic theory. Moreover, we demonstrate how the assumptions on the noise process can be relaxed. Section 4 illustrates the finite sample properties of our approach by means of a Monte Carlo study. Some conclusions and directions for future research are highlighted in Section 5. Finally, we present the proofs in the Appendix.

## 2. Basic notations and definitions

We consider the process $Y$, observed at time points $t_i = i/n$, $i = 0, \ldots, n$. $Y$ is defined on the filtered probability space $(\Omega, \mathcal{F}, (\mathcal{F}_t)_{t \in [0,1]}, P)$ and exhibits a decomposition

$$Y_{i/n} = X_{i/n} + U_i, \tag{2.1}$$

where $X$ is a diffusion process defined by (1.1), and $(U_i)_{0 \leq i \leq n}$ is an i.i.d. noise process with

$$EU_i = 0, \qquad EU_i^2 = \omega^2. \tag{2.2}$$

Further, we assume that $X$ and $U$ are independent.

The core of our approach is the following class of statistics:

$$MBV(Y, r, l)_n = n^{(r+l)/4 - 1/2} \sum_{m=1}^{M} |\bar{Y}_m^{(K)}|^r |\bar{Y}_{m+1}^{(K)}|^l, \qquad r, l \geq 0, \tag{2.3}$$

$$\bar{Y}_m^{(K)} = \frac{1}{n/M - K + 1} \sum_{i=(m-1)n/M}^{mn/M - K} (Y_{(i+K)/n} - Y_{i/n}), \tag{2.4}$$

with

$$K = c_1 n^{1/2}, \qquad M = \frac{n}{c_2 K} = \frac{n^{1/2}}{c_1 c_2} \tag{2.5}$$

for some constants $c_1 > 0$ and $c_2 > 1$ (which will be chosen later).

The intuition behind the quantity $\bar{Y}_m^{(K)}$ can be explained as follows: Since $X$ is a continuous process and $(U_i)_{0 \le i \le n}$ is an i.i.d. process with $EU_i = 0$ we deduce that

$$\frac{1}{n/M - K + 1} \sum_{i=(m-1)n/M}^{mn/M - K} Y_{i/n} = X_{(m-1)/M} + \mathrm{o}_p(1),$$

$$\frac{1}{n/M - K + 1} \sum_{i=(m-1)n/M}^{mn/M - K} Y_{(i+K)/n} = X_{(m-1)/M + K/n} + \mathrm{o}_p(1).$$

This means that the quantity $\bar{Y}_m^{(K)}$ behaves like the increment $X_{(m-1)/M + K/n} - X_{(m-1)/M}$ (although it has a bias that has to be corrected), and consequently it contains information about the volatility process $\sigma$.

The constants $K$ and $M$ control the stochastic order of the term $\bar{Y}_m^{(K)}$. In particular, we have

$$\bar{U}_m^{(K)} = \mathrm{O}_p\left(\sqrt{\frac{1}{n/M - K}}\right) \quad \text{and} \quad \bar{X}_m^{(K)} = \mathrm{O}_p\left(\sqrt{\frac{K}{n}}\right), \tag{2.6}$$

where $U_m^{(K)}$ and $X_m^{(K)}$ are defined analogously to $Y_m^{(K)}$ in (2.4). By (2.5) the stochastic orders of the quantities in (2.6) are balanced, and we obtain

$$\bar{Y}_m^{(K)} = \mathrm{O}_p(n^{-1/4}), \tag{2.7}$$

which explains the normalising factor in (2.3).

More generally, we define the modulated multipower variation by setting

$$MMV(Y, r_1, \ldots, r_k)_n = n^{r_+/4 - 1/2} \sum_{m=1}^{M-k+1} \prod_{j=1}^{k} |\bar{Y}_{m+j-1}^{(K)}|^{r_j},$$

where $k$ is a fixed natural number, $r_j \ge 0$ for all $j$ and $r_+ = r_1 + \cdots + r_k$. This type of construction has been intensively used in a pure Itô diffusion framework (see, e.g., Barndorff-Nielsen, Graversen, Jacod, Podolskij and Shephard [5] or Christensen and Podolskij [11] among others). Later on we will show that the modulated multipower variation for an appropriate choice of $k$ and $r_1, \ldots, r_k$, turns out to be robust to finite activity jumps when estimating arbitrary powers of volatility.

In the sequel we mainly focus on the asymptotic theory of the modulated bipower variation, but we also state the corresponding results for $MMV(Y, r_1, \ldots, r_k)_n$ for the sake of completeness.

# 3. Asymptotic theory

In this section we study the asymptotic behaviour of the class of estimators $MBV(Y, r, l)_n$, $r, l \geq 0$. Before we state the main results of this section we introduce the following notation:

$$\mu_r = E[|z|^r], \qquad z \sim \mathcal{N}(0, 1). \tag{3.1}$$

## 3.1. Consistency

**Theorem 1.** *Assume that* $E|U|^{2(r+l)+\varepsilon} < \infty$ *for some* $\varepsilon > 0$. *If* $M$ *and* $K$ *satisfy* (2.5) *then the convergence in probability*

$$MBV(Y, r, l)_n \xrightarrow{P} MBV(Y, r, l) = \frac{\mu_r \mu_l}{c_1 c_2} \int_0^1 (\nu_1 \sigma_u^2 + \nu_2 \omega^2)^{(r+l)/2} \, du \tag{3.2}$$

*holds. The constants* $\nu_1$ *and* $\nu_2$ *are given by*

$$\nu_1 = \frac{c_1(3c_2 - 4 + (2 - c_2)^3 \vee 0)}{3(c_2 - 1)^2}, \qquad \nu_2 = \frac{2((c_2 - 1) \wedge 1)}{c_1(c_2 - 1)^2}. \tag{3.3}$$

Note that the limit $MBV(Y, r, l)$ in (3.2) depends only on the second moment $\omega^2$ of $U$, and no higher moments are involved.

Next, we present the convergence in probability of the modulated multipower variation $MMV(Y, r_1, \ldots, r_k)_n$.

**Theorem 2.** *Assume that* $E|U|^{2r_+ + \varepsilon} < \infty$ *for some* $\varepsilon > 0$. *If* $M$ *and* $K$ *satisfy* (2.5) *then the convergence in probability*

$$MMV(Y, r_1, \ldots, r_k)_n \xrightarrow{P} MMV(Y, r_1, \ldots, r_k) = \frac{\mu_{r_1} \cdots \mu_{r_k}}{c_1 c_2} \int_0^1 (\nu_1 \sigma_u^2 + \nu_2 \omega^2)^{r_+/2} \, du \tag{3.4}$$

*holds.*

### 3.1.1. Consistent estimates of integrated volatility and integrated quarticity

Theorem 1 shows that $MBV(Y, r, l)_n$ is inconsistent when estimating arbitrary (integrated) powers of volatility. However, when $r + l$ is an even number (this condition is satisfied for the most interesting cases), a slight modification of $MBV(Y, r, l)_n$ turns out to be consistent. Let us illustrate this procedure by providing consistent estimates for the integrated volatility and the integrated quarticity.

As already mentioned in Zhang, Mykland and Ait-Sahalia [25] the statistic

$$\hat{\omega}^2 = \frac{1}{2n} \sum_{i=1}^n |Y_{i/n} - Y_{(i-1)/n}|^2 \tag{3.5}$$

is a consistent estimator of the quantity $\omega^2$ with the convergence rate $n^{-1/2}$. Consequently, we obtain the convergence in probability of the modulated realised volatility

$$MRV(Y)_n := \frac{c_1 c_2 MBV(Y, 2, 0)_n - \nu_2 \hat{\omega}^2}{\nu_1} \xrightarrow{P} \int_0^1 \sigma_u^2 \, \mathrm{d}u \qquad (3.6)$$

as a direct application of Theorem 1 and (3.5).

Now we are in a position to construct a consistent estimator of the integrated quarticity. By (3.6) and Theorem 1 we obtain consistency of the modulated realised quarticity, namely

$$MRQ(Y)_n := \frac{(c_1 c_2 /3) MBV(Y, 4, 0)_n - 2\nu_1 \nu_2 \hat{\omega}^2 MRV(Y)_n - \nu_2^2 (\hat{\omega}^2)^2}{\nu_1^2} \xrightarrow{P} \int_0^1 \sigma_u^4 \, \mathrm{d}u. \qquad (3.7)$$

Note, however, that Theorem 1 provides a whole class of new estimators of the integrated volatility and the integrated quarticity.

### 3.1.2. Robustness to finite activity jumps

As already mentioned in the introduction, one of our main goals is finding consistent estimates of volatility functionals when both microstructure noise and jumps are present. For this purpose we consider the model

$$Z = Y + J, \qquad (3.8)$$

where $Y$ is a noisy diffusion process defined by (2.1) and $J$ denotes a finite activity jump process, that is, $J$ exhibits finitely many jumps on compact intervals. Typical examples of a finite activity jump process are compound Poisson processes.

The next result gives conditions on $r_1, \ldots, r_k$ under which the modulated multipower variation $MMV(Z, r_1, \ldots, r_k)_n$ is robust to finite activity jumps.

**Proposition 3.** *If the assumptions of Theorem 2 are satisfied,* $\max(r_1, \ldots, r_k) < 2$ *and $Z$ is of the form (3.8) then we have*

$$MMV(Z, r_1, \ldots, r_k)_n \xrightarrow{P} MMV(Y, r_1, \ldots, r_k), \qquad (3.9)$$

*where $MMV(Y, r_1, \ldots, r_k)$ is given by (3.4).*

Proposition 3 is shown by the same methods as the corresponding result in the noiseless model (i.e., $U = 0$). We refer to Barndorff-Nielsen, Shephard and Winkel [10] (Proposition 1, page 799) for a detailed proof.

Now we can construct consistent estimates for the integrated volatility and the integrated quarticity, which are robust to noise and finite activity jumps. Since $\hat{\omega}^2$ is robust to jumps, the convergence in probability

$$MBV(Z)_n := \frac{(c_1 c_2 / \mu_1^2) MBV(Z, 1, 1)_n - \nu_2 \hat{\omega}^2}{\nu_1} \xrightarrow{P} \int_0^1 \sigma_u^2 \, \mathrm{d}u \qquad (3.10)$$

holds as a direct consequence of Proposition 3. Similar to the previous subsection, a robust (tripower) estimate of the integrated quarticity is given by

$$
MTQ(Z)_n
$$
$$
:= \frac{(c_1 c_2 / \mu_{2/3}^3) MMV(Z, 4/3, 4/3, 4/3)_n - 2\nu_1 \nu_2 \hat{\omega}^2 MBV(Z)_n - \nu_2^2 (\hat{\omega}^2)^2}{\nu_1^2} \tag{3.11}
$$
$$
\xrightarrow{P} \int_0^1 \sigma_u^4 \, \mathrm{d}u.
$$

**Remark 1.** Recall that the realised volatility $RV$ converges in probability to the integrated volatility plus the sum of squared jumps in the jump-diffusion model. It is interesting to see that the presence of jumps destroys the consistency of the estimator $MRV(Z)_n$, which can be interpreted as an analogue of $RV$. This is explained by the fact that jumps appear with different factors in the statistic $MRV(Z)_n$, according to their positions in the intervals $[\frac{m-1}{M}, \frac{m}{M}]$.

In contrast to our approach, the multiscale estimator of Zhang [24] and the realised kernel estimator of Barndorff-Nielsen, Hansen, Lunde and Shephard [6] converge in probability to the quadratic variation of the jump-diffusion process $X + J$ (in the presence of noise).

Another important object of study is the impact of infinite activity jumps on the modulated bipower (multipower) variation. Such studies can be found in Barndorff-Nielsen, Shephard and Winkel [10], Woerner [23] and in a recent paper of Ait-Sahalia and Jacod [1] for the noiseless framework. We are convinced that similar results hold also for the noisy model, although a more detailed analysis is required.

### 3.1.3. Relaxing the assumptions on the noise process $U$

So far we have assumed that $U$ is an i.i.d. sequence and is independent of the diffusion $X$. Hansen and Lunde [15] have reported that both assumptions are somewhat unrealistic for ultra-high-frequency data. In the following we demonstrate how these conditions can be relaxed.

First, note that the i.i.d. assumption is not essential to guarantee the stochastic order of $\bar{U}_m^{(K)}$ in (2.6). When we assume, for instance, that $U$ is a $q$-dependent sequence, the result of Theorem 1 holds, although higher order autocorrelations of $U$ appear in the limit. In this case we require a stationarity condition on $U$ for the estimation of the autocorrelations and a bias correction of the limit in (3.2).

Further, by using other constants $M$ and $K$ the influence of the noise process $U$ can be made negligible, and independence between $X$ and $U$ is not required. (2.6) implies that, in particular, when we set

$$
K = c_1 n^{1/2+\gamma}, \qquad M = \frac{n}{c_2 K} \tag{3.12}
$$

for some $0 < \gamma < \frac{1}{2}$, the quantity $\bar{X}_m^{(K)}$ driven by the diffusion process dominates the term $\bar{U}_m^{(K)}$. More precisely, the convergence in probability

$$n^{(1-2\gamma)(r+l)/4-(1-2\gamma)/2} \sum_{m=1}^{M} |\bar{Y}_m^{(K)}|^r |\bar{Y}_{m+1}^{(K)}|^l \overset{P}{\longrightarrow} \frac{\mu_r \mu_l \nu_1^{(r+l)/2}}{c_1 c_2} \int_0^1 |\sigma_u|^{r+l}\,\mathrm{d}u \qquad (3.13)$$

holds. The convergence in (3.13) has another useful side effect. It provides consistent estimates for arbitrary integrated powers of volatility. However, since the diffusion term $\bar{X}_m^{(K)}$ dominates the noise term $\bar{U}_m^{(K)}$, the above choice of $K$ and $M$ leads to a slower rate of convergence.

## 3.2. Central limit theorems

In this subsection we present the central limit theorems for a normalised version of $MBV(Y, r, l)_n$. For this purpose we need a structural assumption on the process $\sigma$.

**(V):** The volatility function $\sigma$ satisfies the equation

$$\sigma_t = \sigma_0 + \int_0^t a_s'\,\mathrm{d}s + \int_0^t \sigma_s'\,\mathrm{d}W_s + \int_0^t v_s'\,\mathrm{d}V_s. \qquad (3.14)$$

Here $a'$, $\sigma'$ and $v'$ are adapted cadlag processes, with $a'$ also being predictable and locally bounded, and $V$ is a second Brownian motion independent of $W$.

Condition (V) is a standard assumption that is required for the proof of the central limit theorem for the pure diffusion part $X$ (see, for example, Barndorff-Nielsen, Graversen, Jacod, Podolskij and Shephard [5] or Christensen and Podolskij [11]).

For technical reasons we require a further structural assumption on the noise process $U$. We assume that the filtered probability space $(\Omega, \mathcal{F}, (\mathcal{F}_t)_{t \in [0,1]}, P)$ supports another Brownian motion $B = (B_t)_{t \in [0,1]}$ that is independent of the diffusion process $X$, such that the representation

$$U_i = \sqrt{n}\omega(B_{i/n} - B_{(i-1)/n}) \qquad (3.15)$$

holds.

**Remark 2.** Condition (3.15) ensures that both processes $X$ and $U$ are measurable with respect to the same type of filtration. This assumption enables us to use standard central limit theorems for high frequency observations (see Jacod and Shiryaev [20]). The same assumption has already been used in Gloter and Jacod [13, 14].

The normal distribution of the noise induced by (3.15) is not crucial for our asymptotic theory, and other functions of rescaled increments of $B$ can be considered. Of course, this leads to a slight modification of the central limit theorems presented below.

In the central limit theorems that will be demonstrated below we use the concept of stable convergence of random variables. Let us shortly recall the definition. A sequence of random variables $G_n$ converges stably in law with limit $G$ (throughout this paper we write $G_n \xrightarrow{\mathcal{D}_{st}} G$), defined on an appropriate extension $(\Omega', \mathcal{F}', P')$ of a probability space $(\Omega, \mathcal{F}, P)$, if and only if for any $\mathcal{F}$-measurable and bounded random variable $H$ and any bounded and continuous function $g$ the convergence

$$\lim_{n \to \infty} E[Hg(G_n)] = E[Hg(G)]$$

holds. This is obviously a slightly stronger mode of convergence than convergence in law (see Renyi [22] or Aldous and Eagleson [2] for more details on stable convergence).

Now we present a central limit theorem for the statistic $MBV(Y, r, l)_n$.

**Theorem 4.** *Assume that $U$ is of the form (3.15) and condition* (V) *is satisfied. If $M$ and $K$ satisfy (2.5), and*

1. $r, l \in (1, \infty) \cup \{0\}$ *or*
2. $r$ *or* $l \in (0, 1]$, *and $\sigma_s \neq 0$ for all $s$,*

*then we have*

$$n^{1/4}(MBV(Y, r, l)_n - MBV(Y, r, l)) \xrightarrow{\mathcal{D}_{st}} L(r, l),$$

*where $L(r, l)$ is given by*

$$L(r, l) = \sqrt{\frac{\mu_{2r}\mu_{2l} + 2\mu_r\mu_l\mu_{r+l} - 3\mu_r^2\mu_l^2}{c_1 c_2}} \int_0^1 (\nu_1\sigma_u^2 + \nu_2\omega^2)^{(r+l)/2} \, dW'_u. \qquad (3.16)$$

*Here $W'$ denotes another Brownian motion defined on an extension of the filtered probability space $(\Omega, \mathcal{F}, (\mathcal{F}_t)_{t \in [0,1]}, P)$, which is independent of the $\sigma$-field $\mathcal{F}$.*

Since $\hat{\omega}^2 - \omega^2 = O_p(n^{-1/2})$ we obtain the central limit theorems for the estimates $MRV(Y)_n$ and $MBV(Y)_n$ defined by (3.6) and (3.10), respectively, as a direct consequence of Theorem 4.

**Corollary 1.** *Assume that $U$ is of the form (3.15) and condition* (V) *is satisfied. If $M$ and $K$ satisfy (2.5) then we have*

$$n^{1/4}\left(MRV(Y)_n - \int_0^1 \sigma_u^2 \, du\right) \xrightarrow{\mathcal{D}_{st}} \frac{\sqrt{2c_1 c_2}}{\nu_1} \int_0^1 (\nu_1\sigma_u^2 + \nu_2\omega^2) \, dW'_u, \qquad (3.17)$$

*where $W'$ is another Brownian motion defined on an extension of the filtered probability space $(\Omega, \mathcal{F}, (\mathcal{F}_t)_{t \in [0,1]}, P)$, which is independent of the $\sigma$-field $\mathcal{F}$.*

**Corollary 2.** *Assume that $U$ is of the form (3.15) and condition* (V) *is satisfied. If $M$ and $K$ satisfy (2.5), and $\sigma_s \neq 0$ for all $s$, then we have*

$$n^{1/4}\left(MBV(Y)_n - \int_0^1 \sigma_u^2 \, du\right)$$

$$
\xrightarrow{\mathcal{D}_{st}} \sqrt{\frac{c_1 c_2 (\mu_2^2 + 2\mu_1^2 \mu_2 - 3\mu_1^4)}{\mu_1^4 \nu_1^2}} \int_0^1 (\nu_1 \sigma_u^2 + \nu_2 \omega^2) \, dW_u', \tag{3.18}
$$

*where $W'$ is another Brownian motion defined on an extension of the filtered probability space $(\Omega, \mathcal{F}, (\mathcal{F}_t)_{t \in [0,1]}, P)$, which is independent of the $\sigma$-field $\mathcal{F}$.*

Now let us demonstrate how Corollaries 1 and 2 can be applied in order to obtain confidence intervals for the integrated volatility. Note that the central limit theorem in (3.17) is not feasible yet. Nevertheless, we can easily obtain a feasible version of Corollary 1. Since the Brownian motion $W'$ is independent of the volatility process $\sigma$, the limit defined by (3.17) has a mixed normal distribution with conditional variance

$$
\beta^2 = \frac{2c_1 c_2}{\nu_1^2} \int_0^1 (\nu_1 \sigma_u^2 + \nu_2 \omega^2)^2 \, du.
$$

By an application of Theorem 1 the statistic

$$
\beta_n^2 = \frac{2c_1^2 c_2^2}{3\nu_1^2} MBV(Y, 4, 0)_n
$$

is a consistent estimator of $\beta^2$.

Now we exploit the properties of stable convergence to obtain a standard central limit theorem

$$
\frac{n^{1/4}(MRV(Y)_n - \int_0^1 \sigma_u^2 \, du)}{\beta_n} \xrightarrow{\mathcal{D}} \mathcal{N}(0, 1). \tag{3.19}
$$

From the latter confidence intervals for the integrated volatility can be derived. A feasible version of Corollary 2 can be obtained similarly.

With the above formulae for $\beta^2$ and $\beta_n^2$ in hand we can choose the constants $c_1$ and $c_2$ that minimise the conditional variance. In order to compare our asymptotic variance with the corresponding results of other methods we assume that the volatility process $\sigma$ is constant. In that case the conditional variance $\beta^2$ is minimised by

$$
c_1 = \sqrt{\frac{18}{(c_2 - 1)(4 - c_2)}} \cdot \frac{\omega}{\sigma}, \qquad c_2 = \frac{8}{5}, \tag{3.20}
$$

and is equal to $\frac{256}{3\sqrt{18}} \cdot \sigma^3 \omega \approx 20.11\sigma^3\omega$. Note that the limits in Corollaries 1 and 2 are the same up to a constant. Consequently, the asymptotic conditional variance of $MBV(Y)_n$ is minimised for the same choice of $c_1$ and $c_2$ as above, and is approximately equal to $26.14\sigma^3\omega$, when the volatility function is constant.

As already mentioned in Ait-Sahalia, Mykland and Zhang [25] (see also Gloter and Jacod [13, 14]) the maximum likelihood estimator (when $U$ is normally distributed) converges at the rate $n^{-1/4}$ and has an asymptotic variance $8\sigma^3\omega$, which is a natural lower bound. In contrast to our concept, the family of modified Tukey–Hanning kernel estimators as proposed by Barndorff-Nielsen, Hansen, Lunde and Shephard [6] has an optimal asymptotic variance of about $8.01\sigma^3\omega$. This shows that our estimator is somewhat inefficient in comparison to these kernel-based estimators. A natural direction of future research is to modify our procedure in order to achieve a higher efficiency.

However, the concept of modulated bipower (multipower) variation has been established to provide estimates of arbitrary powers of volatility for the noisy diffusion model, which are additionally robust to finite activity jumps. These are properties which are not captured by the multiscale or realised kernel approach.

For the sake of completeness we state a central limit theorem for the modulated multipower variation $MMV(Y, r_1, \ldots, r_k)_n$.

**Theorem 5.** *Assume that $U$ is of the form (3.15) and condition* (V) *is satisfied. If $M$ and $K$ satisfy (2.5), and*

1. *$r_1, \ldots, r_k \in (1, \infty) \cup \{0\}$ or*
2. *$r_i \in (0, 1]$ for at least one $i$, and $\sigma_s \neq 0$ for all $s$,*

*then we have*

$$n^{1/4}(MMV(Y, r_1, \ldots, r_k)_n - MMV(Y, r_1, \ldots, r_k)) \xrightarrow{\mathcal{D}_{st}} L(r_1, \ldots, r_k),$$

*where $L(r_1, \ldots, r_k)$ is given by*

$$L(r, l) = \sqrt{\frac{A(r_1, \ldots, r_k)}{c_1 c_2}} \int_0^1 (\nu_1 \sigma_u^2 + \nu_2 \omega^2)^{(r+l)/2} \, \mathrm{d}W'_u, \tag{3.21}$$

*with*

$$A(r_1, \ldots, r_k) = \prod_{l=1}^k \mu_{2r_l} - (2k-1) \prod_{l=1}^k \mu_{r_l}^2 + 2 \sum_{j=1}^{k-1} \prod_{l=1}^j \mu_{r_l} \prod_{l=k-j+1}^k \mu_{r_l} \prod_{l=1}^{k-j} \mu_{r_l + r_{l+j}}.$$

*Here $W'$ denotes another Brownian motion defined on an extension of the filtered probability space $(\Omega, \mathcal{F}, (\mathcal{F}_t)_{t \in [0,1]}, P)$, which is independent of the $\sigma$-field $\mathcal{F}$.*

Note that the constant $A(r_1, \ldots, r_k)$ also appears in the central limit theorem for multipower variation in a pure diffusion framework (see Barndorff-Nielsen, Graversen, Jacod, Podolskij and Shephard [5]).

## 4. Simulation study

In this section, we inspect the finite sample properties of various proposed estimators for both integrated volatility and quarticity through Monte Carlo experiments. Moreover, we compare our estimators' behaviour with the properties of the corresponding kernel-based estimators from Barndorff-Nielsen, Hansen, Lunde and Shephard [6]. To this end, we choose the same stochastic volatility model as in their work, namely

$$
\begin{aligned}
\mathrm{d}X_t &= \mu\,\mathrm{d}t + \sigma_t\,\mathrm{d}W_t, & \sigma_t &= \exp(\beta_0 + \beta_1\tau_t), \\
\mathrm{d}\tau_t &= \alpha\tau_t\,\mathrm{d}t + \mathrm{d}B_t, & \mathrm{corr}(\mathrm{d}W_t, \mathrm{d}B_t) &= \rho,
\end{aligned}
\tag{4.1}
$$

with $\mu = 0.03, \beta_0 = 0.3125, \beta_1 = 0.125, \alpha = -0.025$ and $\rho = -0.3$. $(U_i)_{0 \leq i \leq n}$ is assumed to be i.i.d. normally distributed with variance $\omega^2$.

## 4.1. Simulation design

We create $20\,000$ repetitions of the system in equation (4.1), for which we use an Euler approximation and different values of $n$. Whenever we have to estimate $\omega^2$, we choose $\hat{\omega}^2$ as defined in (3.5).

Since we state propositions for a whole class of estimators, we do not focus on one special estimator. To be precise, we investigate the finite sample properties in three different situations.

First, we study the performance of $MRV(Y)_n$ as an estimator for the integrated volatility and compare it with the corresponding kernel-based statistic from Barndorff-Nielsen, Hansen, Lunde and Shephard [6], using the modified Tukey–Hanning kernel with $p = 2$. We denote this estimator by $KB(Y)_n$. In Table 1 we present the Monte Carlo results for both mean and variance of the two statistics for various choices of $n$ and $\omega^2 = 0.01, 0.001$, which is a reasonable choice, since $IV$ is about 2 in model (4.1). Figure 1 gives histograms both of the standardised statistic in (3.19) and of the statistic

$$
\frac{n^{1/4}(\log(MRV(Y)_n) - \log(\int_0^1 \sigma_u^2\,\mathrm{d}u))}{\beta_n / MRV(Y)_n} \xrightarrow{\mathcal{D}} \mathcal{N}(0,1),
\tag{4.2}
$$

which is obtained via an application of the delta method. Both statistics converge weakly to a standard normal distribution. In this case, we choose two different values of $n$, namely $n = 1024$ and $n = 16\,384$.

Second, we analyse the performance of the estimation of the integrated volatility in the presence of finite activity jumps. In this case we use the bipower estimator $MBV(Z)_n$, which is robust to jumps. We present the Monte Carlo results for this estimator in Table 2.

At last, we analyse how well $MRQ(Y)_n$ works as an estimator for the integrated quarticity in contrast to the proposed bipower variation estimator in Barndorff-Nielsen, Hansen, Lunde and Shephard [6], which we call $BP(Y)_n$. Note that $BP(Y)_n$ has a convergence rate of $n^{-1/6}$, which is obviously slower than the convergence rate of our estimator $MRQ(Y)_n$. Table 3 shows the results in the quite simple setting

$$
\mathrm{d}X_t = \mu\,\mathrm{d}t + \mathrm{d}W_t
\tag{4.3}
$$

with $\mu = 0.03$ as above.

As mentioned before, the asymptotic (conditional) variance of the estimators $MRV(Y)_n$ and $MBV(Y)_n$ can be minimised for an appropriate choice of $c_1$ and $c_2$, if the volatility function is constant and the drift function is zero. However, even in model (4.1) and (4.3) the choice of $c_1$ and $c_2$ as in (3.20) and with $\sigma$ replaced by $IV$ may give an idea of a reasonable size for the two constants. Using a first estimate for $IV$ we decided to choose $c_1 = 0.25$ for $\omega^2 = 0.01$ and $c_1 = 0.125$ for $\omega^2 = 0.001$, whereas $c_2 = 2$. Since the computation of the optimal values of $c_1$ and $c_2$ for the estimation of $IQ$ involves the solution of polynomial equations with higher degrees than two, we have dispensed with this analysis and set $c_1 = 1$ and $c_2 = 1.6$, both for $\omega^2 = 0.01$ and $\omega^2 = 0.001$. To produce the process $J$ we allocate one jump in the interval $[0,1]$. The arrival time of this jump is considered to be uniformly distributed, whereas the jump size is $\mathcal{N}(0, h^2)$ distributed with $h = 0.1, 0.25$.

Note finally that it might be convenient not to plug in $\nu_1$ to compute $MRV(Y)_n$ and the estimator for the conditional variance $\beta_n^2$, but to use

$$\nu_1^{(n)} = \nu_1 + \frac{(3 - c_2) \wedge 1/(c_2 - 1)}{(c_2 - 1)\sqrt{n}} + \mathrm{O}\left(\frac{1}{n}\right),$$

which is a better approximation to the second moment of $n^{1/4}\bar{W}_m^{(K)}$ than $\nu_1$. Similarly, one could use a finite sample analogue $\nu_2^{(n)}$ for $\nu_2$ as well. It has turned out that the performance of our estimators is indeed sensitive to the choice of $\nu_1$ (which is why we used $\nu_1^{(n)}$), whereas it is almost unaffected by the transition from $\nu_2$ to $\nu_2^{(n)}$ (which is why we used $\nu_2$).

## 4.2. Results

Since our aim is mainly to give an idea of how well the different estimators work, we content ourselves with computing the estimated mean and variance of the bias-corrected statistics. Except for $MRV(Y)_n$ we therefore do not evaluate the accuracy of the stated central limit theorems.

Table 1 shows that $MRV(Y)_n$ works quite well as an estimator of the integrated volatility in the noisy diffusion setting, since both bias and variance are rather small, at least for sample sizes larger than $n = 1024$. For large values of $n$ and $\omega^2 = 0.01$ it provides even better finite sample properties than $KB(Y)_n$, whereas the kernel-based estimator improves a lot when the variance of the noise terms becomes smaller. Nevertheless, $MRV(Y)_n$ is a serious alternative to the kernel-based estimator, especially for large values of $\omega^2$.

Figure 1 indicates that the behaviour of each standardised statistic depends heavily on the actual size of $\omega^2$, especially for a small sample size $n$. In fact, even for $n = 1024$ we see a reasonable approximation of the standard normal density when the variance of the noise variables is large, whereas for $\omega^2 = 0.001$ the histogram exhibits a significant shift to the right. One might suggest that these effects are caused by a large variance of the estimator of the integrated quarticity. For a large value of $n$ all statistics work pretty well; however, it is remarkable that the transition to the log-transformed statistic provides an improvement, at least for a large choice of $\omega^2$.

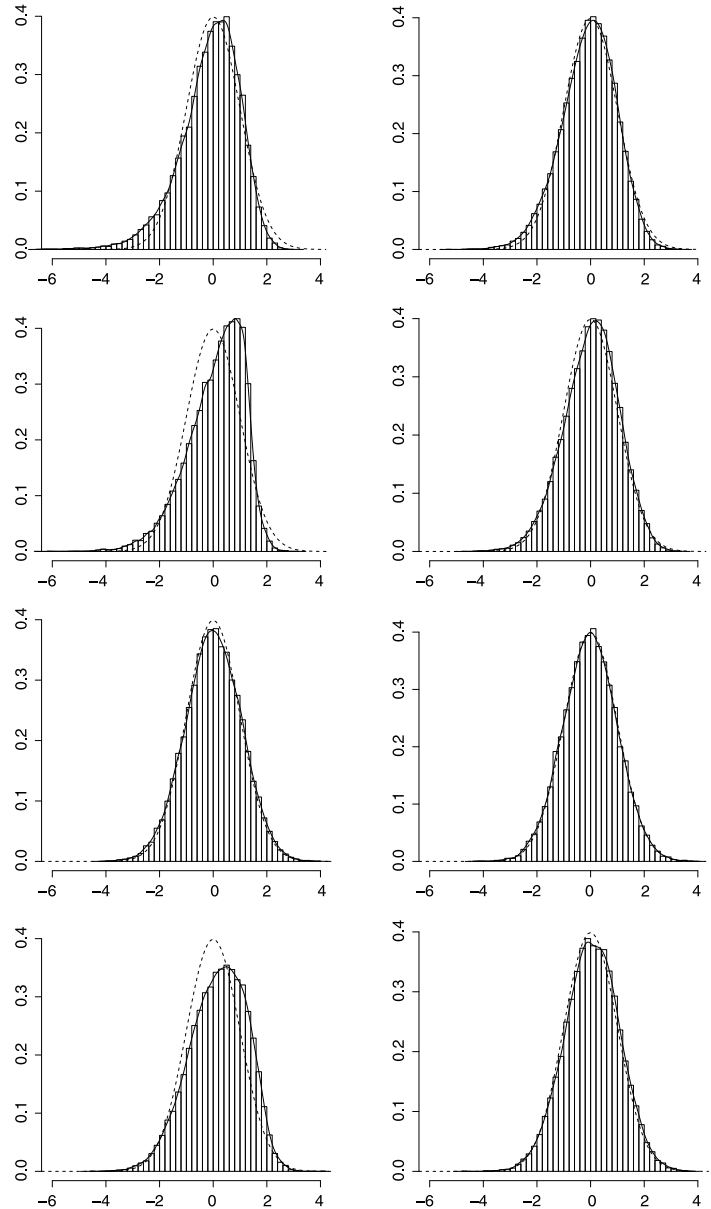

**Figure 1.** Various histograms of the statistics defined in (3.19) and (4.2). In each line, the first data set was computed with $n = 1024$, whereas for the second one we used $n = 16\,384$. The first four histograms illustrate the case (3.19) with $\omega^2 = 0.01$ (in the first line) and $\omega^2 = 0.001$ (in the second line). The latter four graphics show the finite sample properties for the weak convergence in (4.2), in the same order as above. For comparison, the dashed line shows the graph of the standard normal density and the solid line gives a standard kernel density estimate.

**Table 1.** The Monte Carlo results for mean and variance of both $MRV(Y)_n - \int_0^1 \sigma_u^2 \, du$ and $KB(Y)_n - \int_0^1 \sigma_u^2 \, du$ for various values of $n$ and $\omega^2$. The data are generated from the model (4.1)

| $n$ | $\omega^2 = 0.01$ | | $\omega^2 = 0.001$ | |
|---|---|---|---|---|
| | Mean | Variance | Mean | Variance |
| $MRV(Y)_n$ | | | | |
| 256 | 0.1363 | 0.63 | 0.5245 | 1.782 |
| 1024 | 0.0433 | 0.219 | 0.1717 | 0.269 |
| 4096 | 0.0113 | 0.102 | 0.0478 | 0.055 |
| 9216 | 0.0045 | 0.064 | 0.0243 | 0.031 |
| 16 384 | 0.0059 | 0.05 | 0.0129 | 0.021 |
| 25 600 | 0.004 | 0.039 | 0.0094 | 0.017 |
| $KB(Y)_n$ | | | | |
| 256 | −0.022 | 0.228 | −0.0289 | 0.143 |
| 1024 | 0.0074 | 0.091 | −0.0075 | 0.042 |
| 4096 | 0.0195 | 0.046 | 0.0004 | 0.015 |
| 9216 | 0.0203 | 0.038 | 0.001 | 0.009 |
| 16 384 | 0.0201 | 0.04 | 0.001 | 0.007 |
| 25 600 | 0.0178 | 0.046 | 0.0013 | 0.005 |

**Table 2.** Mean and variance of $MBV(Z)_n - \int_0^1 \sigma_u^2 \, du$ in the presence of jumps. We choose the sample frequency as before and analyse the finite sample properties for different values of $\omega^2$ and $h$, where $h$ denotes the variance of the jump size

| $n$ | $\omega^2 = 0.01,\ h = 0.25$ | | $\omega^2 = 0.001,\ h = 0.25$ | | $\omega^2 = 0.001,\ h = 0.1$ | |
|---|---|---|---|---|---|---|
| | Mean | Variance | Mean | Variance | Mean | Variance |
| 256 | 0.0582 | 0.614 | −0.0839 | 0.332 | −0.1224 | 0.29 |
| 1024 | 0.0835 | 0.295 | 0.0274 | 0.133 | −0.102 | 0.112 |
| 4096 | 0.0707 | 0.15 | 0.0466 | 0.063 | 0.0184 | 0.056 |
| 9216 | 0.0642 | 0.102 | 0.0461 | 0.043 | 0.0107 | 0.038 |
| 16 384 | 0.0599 | 0.076 | 0.044 | 0.032 | 0.025 | 0.028 |
| 25 600 | 0.0566 | 0.059 | 0.0415 | 0.025 | 0.0181 | 0.023 |

From Table 2 we conclude that in the noisy jump-diffusion framework the proposed bipower estimator $MBV(Z)_n$ works quite well. Both bias and the variance of $MBV(Z)_n$ are rather small, even for moderate values of $n$.

Finally, we see from Table 3 that $MRQ(Y)_n$ is on average closer to the true integrated quarticity than $BP(Y)_n$, whereas the variance of $BP(Y)_n$ is smaller than that of $MRQ(Y)_n$, even though $BP(Y)_n$ has a slower rate of convergence. However, we are convinced that the efficiency of $MRQ(Y)_n$ can be improved massively by choosing the constants $c_1$ and $c_2$ optimally.

## 5. Conclusions and directions for future research

This paper highlights the potential of the modulated bipower approach, and we are convinced that many unsolved problems in a noisy (jump-)diffusion framework can be tackled by our methods. Let us mention some most important directions for future research. First, we intend to modify our approach by subsampling the statistic $MBV(Y, r, l)_n$ to obtain more efficient estimators of the integrated volatility and the integrated quarticity. A first step in this direction has been made in a recent paper by Jacod, Li, Mykland, Podolskij and Vetter [18] who proposed a subsampled version of $MBV(Y, 2, 0)_n$ to estimate $IV$ in the presence of a more general noise process. Second, we plan to derive a multivariate version of the current approach. This can be used to estimate the quadratic covariation, which is a key concept in econometrics, in the presence of noise. Finally, an interesting and very important modification of this problem is the estimation of the quadratic covariation for non-synchronously observed data in the presence of noise (see Hayashi and Yoshida [16] for more details in a pure diffusion case).

**Table 3.** The finite sample properties of $MRQ(Y)_n - \int_0^1 \sigma_u^4 \, du$ and $BP(Y)_n - \int_0^1 \sigma_u^4 \, du$ in model (4.3). Both sample frequency and noise are the same as in Table 1

| $n$ | $\omega^2 = 0.01$ | | $\omega^2 = 0.001$ | |
|---|---|---|---|---|
| | Mean | Variance | Mean | Variance |
| $MRQ(Y)_n$ | | | | |
| 256 | 0.0745 | 1.348 | 0.0686 | 1.274 |
| 1024 | 0.0128 | 0.587 | 0.0121 | 0.557 |
| 4096 | 0.0135 | 0.306 | 0.0013 | 0.278 |
| 9216 | 0.0113 | 0.203 | 0.015 | 0.184 |
| 16 384 | 0.0159 | 0.152 | 0.0155 | 0.14 |
| 25 600 | 0.0088 | 0.117 | 0.0077 | 0.108 |
| $BP(Y)_n$ | | | | |
| 256 | −0.2517 | 0.304 | −0.2803 | 0.274 |
| 1024 | −0.1811 | 0.186 | −0.1434 | 0.169 |
| 4096 | −0.0312 | 0.108 | −0.0745 | 0.095 |
| 9216 | −0.0089 | 0.077 | −0.04 | 0.065 |
| 16 384 | 0.0078 | 0.059 | −0.0287 | 0.048 |
| 25 600 | 0.0148 | 0.047 | −0.0206 | 0.039 |

# Appendix

In the following we assume without loss of generality that $a$, $\sigma$, $a'$, $\sigma'$ and $v'$ are bounded (for details see Barndorff-Nielsen, Graversen, Jacod, Podolskij and Shephard [5]). Moreover, the constants that appear in the proofs are all denoted by $C$ (even if they depend on the actual choice of the powers $l$ and $r$).

Note first that the numbers $\nu_1$ and $\nu_2$ as given in (3.3) are the limits of the variances of the random variables $n^{1/4}\bar{W}_m^{(K)}$ and $n^{1/4}\bar{U}_m^{(K)}$, respectively. Since a standard calcuation shows that the true variances $\nu_i^{(n)}$ are of the form

$$\nu_i^{(n)} = \nu_i + \mathrm{O}(n^{-1/2}),$$

one may conclude from the mean value theorem that replacing $\nu_i^{(n)}$ by $\nu_1$ does affect neither the consistency results nor the central limit theorem. Thus, whenever we refer to the variances of those random quantities, we use $\nu_1$ and $\nu_2$ without further notice.

Before we start with the proofs of main results, we introduce more notation and state some simple lemmata. We consider the quantities

$$\beta_m^n = n^{1/4}(\sigma_{(m-1)/M}\bar{W}_m^{(K)} + \bar{U}_m^{(K)}), \qquad \beta_m'^n = n^{1/4}(\sigma_{(m-1)/M}\bar{W}_{m+1}^{(K)} + \bar{U}_{m+1}^{(K)}), \quad (6.1)$$

which approximate $\bar{Y}_m^{(K)}$ and $\bar{Y}_{m+1}^{(K)}$, respectively, by using the associated increments of the underlying Brownian motion $W$. We further define

$$\xi_m^n = n^{1/4}\bar{Y}_m^{(K)} - \beta_m^n, \qquad \xi_m'^n = n^{1/4}\bar{Y}_{m+1}^{(K)} - \beta_m'^n \qquad (6.2)$$

as the differences between the true quantities and their approximations. We further set $f(x) := |x|^r$ and $g(x) := |x|^l$. In the next lemma we study the stochastic order of the terms $\beta_m^n$ and $\xi_m^n$.

**Lemma 1.** *We have*

$$E[|\xi_m^n|^q] + E[|\xi_m'^n|^q] + E[|n^{1/4}\bar{X}_m^{(K)}|^q] < C \qquad (6.3)$$

*for any $q > 0$, and*

$$E[|\beta_m^n|^q] + E[|\beta_m'^n|^q] + E[|n^{1/4}\bar{Y}_m^{(K)}|^q] < C \qquad (6.4)$$

*for any $0 < q < 2(r+l) + \varepsilon$ with $\varepsilon$ as stated in Theorem 1. Both results hold uniformly in $m$.*

**Proof.** We begin with the proof of (6.3). In the case $q \geq 1$, we use Hölder's inequality to obtain

$$E[|\xi_m^n|^q]$$

$$= E\left[\left|\frac{n^{1/4}}{n/M - K + 1}\sum_{i=n(m-1)/M}^{nm/M - K}(X_{(i+K)/n} - X_{i/n}) - \sigma_{(m-1)/M}(W_{(i+K)/n} - W_{i/n})\right|^q\right]$$

$$\leq \frac{1}{n/M - K + 1}$$

$$\times \sum_{i=n(m-1)/M}^{nm/M - K} E\left[\left|n^{1/4}((X_{(i+K)/n} - X_{i/n}) - \sigma_{(m-1)/M}(W_{(i+K)/n} - W_{i/n}))\right|^q\right]$$

$$= \frac{1}{n/M - K + 1}$$

$$\times \sum_{i=n(m-1)/M}^{nm/M - K} E\left[\left|n^{1/4}\left(\int_{i/n}^{(i+K)/n} a_s\,\mathrm{d}s + \int_{i/n}^{(i+K)/n}(\sigma_s - \sigma_{(m-1)/M})\,\mathrm{d}W_s\right)\right|^q\right].$$

Thus the property follows from the boundedness of the functions $a$ and $\sigma$, and a use of Burkholder's inequality. For $q < 1$, Jensen's inequality yields

$$E[|\xi_m^n|^q] \leq E[|\xi_m^n|]^q,$$

and we have (6.3) just as above. The corresponding assertions for $\xi_m'^n$ and $n^{1/4}\bar{X}_m^{(K)}$ can be shown analogously.

Now let us prove (6.4). In the same way as before we have

$$E[|n^{1/4}\bar{Y}_m^{(K)}|^q] \leq C(E[|n^{1/4}\bar{U}_m^{(K)}|^q] + E[|n^{1/4}\bar{X}_m^{(K)}|^q])$$

for any $q \geq 0$. Moreover, it can be shown that $n^{1/4}\bar{U}_m^{(K)}$ can be rewritten as a weighted sum of independent random variables, for which the convergence in distribution

$$n^{1/4}\bar{U}_m^{(K)} \xrightarrow{\mathcal{D}} \mathcal{N}(0, \nu_2\omega^2)$$

holds. Using the continuity theorem and the moment assumption for each $0 < q < 2(r + l) + \varepsilon$ we obtain by uniform integrability of $|n^{1/4}\bar{U}_m^{(K)}|^q$ that $E[|n^{1/4}\bar{U}_m^{(K)}|^q]$ is bounded. This proves (6.4) for $n^{1/4}\bar{Y}_m^{(K)}$. The corresponding result for the quantities $\beta_m^n$ and $\beta_m'^n$ can be shown analogously. □

The next lemma will be used later to prove that the error due to the approximation of $\bar{Y}_m^{(K)}$ using $\beta_m^n$ is small compared to the rate of convergence. For a more general setting and a proof see Lemma 5.4 in Barndorff-Nielsen, Graversen, Jacod, Podolskij and Shephard [5].

**Lemma 2.** *If*

$$Z_m^n := 1 + |\mu_m^n| + |\mu_m'^n| + |\mu_m''^n|$$

*satisfies $E[|Z_m^n|^q] < C$ for all $0 < q < 2(r+l) + \varepsilon$ and if further*

$$\frac{1}{M} \sum_{m=1}^{M} E[|\mu_m'^n - \mu_m''^n|^2] \to 0 \tag{6.5}$$

*holds, then we have*

$$\frac{1}{M} \sum_{m=1}^{M} E[f^2(\mu_m^n)(g(\mu_m'^n) - g(\mu_m''^n))^2] \to 0.$$

**Proof of Theorem 1.** We introduce the quantities

$$MBV^n := \sum_{m=1}^{M} \eta_m^n \quad \text{and} \quad MBV'^n := \sum_{m=1}^{M} \eta_m'^n,$$

where $\eta_m^n$ and $\eta_m'^n$ are defined by

$$\eta_m^n := \frac{n^{(r+l)/4}}{c_1 c_2} E[|\bar{Y}_m^{(K)}|^r |\bar{Y}_{m+1}^{(K)}|^l | \mathcal{F}_{(m-1)/M}], \qquad \eta_m'^n := \frac{\mu_r \mu_l}{c_1 c_2} (\nu_1 \sigma_{(m-1)/M}^2 + \nu_2 \omega^2)^{(r+l)/2}.$$

Riemann integrability yields

$$\frac{1}{M} MBV'^n \overset{P}{\longrightarrow} MBV(Y, r, l),$$

so we are forced to prove

$$MBV(Y, r, l)_n - \frac{1}{M} MBV^n \overset{P}{\longrightarrow} 0 \tag{6.6}$$

and

$$\frac{1}{M}(MBV^n - MBV'^n) \overset{P}{\longrightarrow} 0 \tag{6.7}$$

in two steps.

Considering the first step we recall the identity $\sqrt{n} = c_1 c_2 M$ and obtain therefore

$$MBV(Y, r, l)_n - \frac{1}{M} MBV^n = \sum_{m=1}^{M} (\gamma_m - E[\gamma_m | \mathcal{F}_{(m-1)/M}]),$$

where $\gamma_m$ is given by

$$\gamma_m = n^{(r+l)/4 - 1/2} |\bar{Y}_m^{(K)}|^r |\bar{Y}_{m+1}^{(K)}|^l.$$

Using Lenglart's inequality and the $\mathcal{F}_{(m+1)/M}$-measurablity of $\gamma_m$ (for details see Lemma 5.2 in Barndorff-Nielsen, Graversen, Jacod, Podolskij and Shephard [5]) we find that the stochastic convergence stated in (6.6) follows from

$$\sum_{m=1}^{M} E[|\gamma_m|^2 | \mathcal{F}_{(m-1)/M}] \xrightarrow{P} 0.$$

Hölder's inequality and Lemma 1 yield $E[|\gamma_m|^2 | \mathcal{F}_{(m-1)/M}] \leq Cn^{-1}$, thus (6.6) follows.

To prove the assertion in (6.7) recall that $f(x) = |x|^r$ and $g(x) = |x|^l$ and observe that the continuity theorem implies

$$E[n^{(r+l)/4} f(\sigma_{(m-1)/M} \bar{W}_m^{(K)} + \bar{U}_m^{(K)}) g(\sigma_{(m-1)/M} \bar{W}_{m+1}^{(K)} + \bar{U}_{m+1}^{(K)}) | \mathcal{F}_{(m-1)/M}]$$
$$= \mu_r \mu_l (\nu_1 \sigma_{(m-1)/M}^2 + \nu_2 \omega^2)^{(r+l)/2} + o_p(1),$$

uniformly in $m$. Thus

$$\frac{1}{M}(MBV^n - MBV'^n) = \frac{1}{M} \sum_{m=1}^{M} E[\zeta_m^n | \mathcal{F}_{(m-1)/M}] + o_p(1)$$

with

$$\zeta_m^n = \frac{n^{(r+l)/4}}{c_1 c_2}(f(\bar{Y}_m^{(K)}) g(\bar{Y}_{m+1}^{(K)}) - f(\sigma_{(m-1)/M} \bar{W}_m^{(K)} + \bar{U}_m^{(K)}) g(\sigma_{(m-1)/M} \bar{W}_{m+1}^{(K)} + \bar{U}_{m+1}^{(K)})).$$

To obtain the desired result it suffices to show

$$\frac{1}{M} \sum_{m=1}^{M} E[|\zeta_m^n|] \to 0,$$

which follows from

$$\frac{1}{M} \sum_{m=1}^{M} E[|\zeta_m^n|^2] \to 0 \tag{6.8}$$

using the Cauchy–Schwarz inequality. In a first step we obtain for some constant $C > 0$

$$|\zeta_m^n|^2 = \frac{1}{c_1^2 c_2^2}(f(\xi_m^n + \beta_m^n) g(\xi_{m+1}^n + \beta_{m+1}^n) - f(\beta_m^n) g(\beta_m'^n))^2$$

$$\leq C(g^2(\xi_{m+1}^n + \beta_{m+1}^n)(f(\xi_m^n + \beta_m^n) - f(\beta_m^n))^2$$
$$+ f^2(\beta_m^n)(g(\xi_{m+1}^n + \beta_{m+1}^n) - g(\beta_{m+1}'^n))^2 + f^2(\beta_m^n)(g(\beta_{m+1}'^n) - g(\beta_m'^n))^2), \tag{6.9}$$

where the quantities $\beta_m^n$ and $\xi_m^n$ are defined by (6.1) and (6.2), respectively. Since we have shown in (6.3) and (6.4) that the conditions on the boundedness of $Z_m^n$ for an application of Lemma 2 are fulfilled, it suffices to prove

$$\frac{1}{M}\sum_{m=1}^{M}E[|\xi_m^n|^2 + |\xi_{m+1}^n|^2 + |\beta_{m+1}^n - \beta_m'^n|^2] \to 0 \tag{6.10}$$

to obtain the assertion.

For the first term in (6.10) we have

$$E[|\xi_m^n|^2] \leq \frac{1}{n/M - K + 1}$$

$$\times \sum_{i=n(m-1)/M}^{n/M-K} E[|n^{1/4}((X_{(i+K)/n} - X_{i/n}) - \sigma_{(m-1)/M}(W_{(i+K)/n} - W_{i/n}))|^2]$$

as in the proof of (6.3). Using (2.5) and

$$(X_{(i+K)/n} - X_{i/n}) - \sigma_{(m-1)/M}(W_{(i+K)/n} - W_{i/n})$$

$$= \int_{i/n}^{(i+K)/n} a_s \,\mathrm{d}s + \int_{i/n}^{(i+K)/n} (\sigma_s - \sigma_{(m-1)/M}) \,\mathrm{d}W_s$$

we obtain

$$E[|n^{1/4}((X_{(i+K)/n} - X_{i/n}) - \sigma_{(m-1)/M}(W_{(i+K)/n} - W_{i/n}))|^2]$$

$$\leq C\left(n^{-1/2} + n^{1/2}E\left[\int_{i/n}^{(i+K)/n} (\sigma_s - \sigma_{(m-1)/M})^2 \,\mathrm{d}s\right]\right)$$

$$\leq C\left(n^{1/2}E\left[\int_{(m-1)/M}^{m/M} (\sigma_s - \sigma_{(m-1)/M})^2 \,\mathrm{d}s\right]\right) + \mathrm{o}(1).$$

Consequently,

$$\frac{1}{M}\sum_{m=1}^{M}E[|\xi_m^n|^2] \leq C\sum_{m=1}^{M}E\left[\int_{(m-1)/M}^{m/M} (\sigma_s - \sigma_{(m-1)/M})^2 \,\mathrm{d}s\right] + \mathrm{o}(1)$$

$$= C\sum_{m=1}^{M}E\left[\int_{(m-1)/M}^{m/M} (\sigma_s - \sigma_{\lfloor Ms\rfloor/M})^2 \,\mathrm{d}s\right] + \mathrm{o}(1)$$

$$= C\int_0^1 E[(\sigma_s - \sigma_{\lfloor Ms\rfloor/M})^2] \,\mathrm{d}s + \mathrm{o}(1)$$

follows. Since $\sigma$ is bounded and cadlag, Lebesgue's theorem yields

$$\frac{1}{M}\sum_{m=1}^{M}E[|\xi_m^n|^2] \to 0$$

and as well for the second term in (6.10). We further have

$$\beta_{m+1}^n - \beta_m'^n = n^{1/4}(\sigma_{m/M} - \sigma_{(m-1)/M})\bar{W}_{m+1}^{(K)}.$$

Since $\bar{W}_{m+1}^{(K)}$ is independent of $\sigma_t$ for any $t \leq \frac{m}{M}$ we obtain

$$\frac{1}{M}\sum_{m=1}^{M} E[|\beta_{m+1}^n - \beta_m'^n|^2] \leq \frac{C}{M}\sum_{m=1}^{M} E[|\sigma_{m/M} - \sigma_{(m-1)/M}|^2]$$

$$\leq \frac{C}{M}\sum_{m=1}^{M} E[|\sigma_{m/M} - \sigma_s|^2 + |\sigma_s - \sigma_{(m-1)/M}|^2].$$

The assertion therefore follows with the same arguments as above. That completes the proof of (6.7). □

**Proof of Theorem 2.** Theorem 2 can be proven by the same methods as Theorem 1. □

**Proof of Theorem 4.** Here we mainly use the same techniques as presented in Barndorff-Nielsen, Graversen, Jacod, Podolskij and Shephard [5] or Christensen and Podolskij [11]. We will state the proof of the key steps and refer to the articles quoted above for the details.

We define the quantity

$$L_n(r,l) = n^{-1/4}\sum_{m=1}^{M}(f(\beta_m^n)g(\beta_m'^n) - E[f(\beta_m^n)g(\beta_m'^n)|\mathcal{F}_{(m-1)/M}]), \qquad (6.11)$$

where the terms $\beta_m^n$ and $\beta_m'^n$ are given by (6.1), and $f(x) = |x|^r$, $g(x) = |x|^l$. In the next lemma we state the central limit theorem for $L_n(r,l)$.

**Lemma 3.** *We have*

$$L_n(r,l) \xrightarrow{\mathcal{D}_{st}} L(r,l),$$

*where $L(r,l)$ is defined in Theorem 4.*

**Proof.** First, note that

$$L_n(r,l) = \sum_{m=2}^{M+1} \theta_m^n + o_p(1),$$

where $\theta_m^n$ is given by

$$\theta_m^n = n^{-1/4}(f(\beta_{m-1}^n)(g(\beta_{m-1}'^n) - \mu_l(\nu_1\sigma_{(m-2)/M}^2 + \nu_2\omega^2)^{l/2})$$

$$+ \mu_l(\nu_1\sigma_{(m-1)/M}^2 + \nu_2\omega^2)^{l/2}(f(\beta_m^n) - \mu_r(\nu_1\sigma_{(m-1)/M}^2 + \nu_2\omega^2)^{r/2})).$$

We have that

$$E[\theta_m^n|\mathcal{F}_{(m-1)/M}] = 0,$$

and

$$\sum_{m=2}^{M+1} E[|\theta_m^n|^2|\mathcal{F}_{(m-1)/M}] \xrightarrow{P} \frac{\mu_{2r}\mu_{2l} + 2\mu_r\mu_l\mu_{r+l} - 3\mu_r^2\mu_l^2}{c_1 c_2} \int_0^1 (\nu_1\sigma_u^2 + \nu_2\omega^2)^{r+l} \, \mathrm{d}u.$$

Next, let $Z = W$ or $B$. Since $\theta_m^n$ is an even functional in $W$ and $B$, and $(W, B) \overset{\mathcal{D}}{=} -(W, B)$, we obtain the identity

$$E[\theta_m^n(Z_{m/M} - Z_{(m-1)/M})|\mathcal{F}_{(m-1)/M}] = 0.$$

Finally, let $N = (N_t)_{t \in [0,1]}$ be a bounded martingale on $(\Omega, \mathcal{F}, (\mathcal{F}_t)_{t \in [0,1]}, P)$, which is orthogonal to $W$ and $B$ (i.e., with quadratic covariation $[W, N]_t = [B, N]_t = 0$ almost surely). By the arguments of Barndorff-Nielsen, Graversen, Jacod, Podolskij and Shephard [5] we have

$$E[\theta_m^n(N_{m/M} - N_{(m-1)/M})|\mathcal{F}_{(m-1)/M}] = 0.$$

Now the stable convergence in Lemma 3 follows by Theorem IX 7.28 in Jacod and Shiryaev [20]. □

Now we are left to prove the convergence

$$n^{1/4}(MBV(Y,r,l)_n - MBV(Y,r,l)) - L_n(r,l) \xrightarrow{P} 0. \tag{6.12}$$

Obviously, the convergence in (6.12) is equivalent to

$$\sum_{m=1}^{M} \vartheta_m^n \xrightarrow{P} 0, \tag{6.13}$$

$$\sum_{m=1}^{M} \vartheta_m'^n \xrightarrow{P} 0, \tag{6.14}$$

with $\vartheta_m^n, \vartheta_m'^n$ defined by

$$\vartheta_m^n = n^{-1/4}[f(n^{1/4}\bar{Y}_m^{(K)})g(n^{1/4}\bar{Y}_{m+1}^{(K)}) - f(\beta_m^n)g(\beta_m'^n)],$$

$$\vartheta_m'^n = n^{1/4} \int_{(m-1)/M}^{m/M} ((\nu_1\sigma_u^2 + \nu_2\omega^2)^{(r+l)/2} - (\nu_1\sigma_{(m-1)/M}^2 + \nu_2\omega^2)^{(r+l)/2}) \, \mathrm{d}u.$$

The convergence in (6.14) has been shown in Barndorff-Nielsen, Graversen, Jacod, Podolskij and Shephard [5], and so we concentrate on proving (6.13). For the same reason as in the proof of (6.6) (and from a similar argument as in the proof of (6.8)) this result follows from

$$\sum_{m=1}^{M} E[\vartheta_m^n | \mathcal{F}_{(m-1)/M}] \xrightarrow{P} 0. \tag{6.15}$$

Observe that

$$\vartheta_m^n = n^{-1/4} f(n^{1/4} \bar{Y}_m^{(K)})(g(n^{1/4} \bar{Y}_{m+1}^{(K)}) - g(\beta_m'^n)) + n^{-1/4} g(\beta_m'^n)(f(n^{1/4} \bar{Y}_m^{(K)}) - f(\beta_m^n)).$$

Now we obtain

$$\sum_{m=1}^{M} E[\vartheta_m^n | \mathcal{F}_{(m-1)/M}] = \sum_{m=1}^{M} E[\vartheta_m^n(1) + \vartheta_m^n(2) | \mathcal{F}_{(m-1)/M}] + o_p(1), \tag{6.16}$$

with $\vartheta_m^n(1)$, $\vartheta_m^n(2)$ defined by

$$\vartheta_m^n(1) = n^{-1/4} \nabla g(\beta_m'^n) f(n^{1/4} \bar{Y}_m^{(K)}) \xi_m'^n,$$
$$\vartheta_m^n(2) = n^{-1/4} \nabla f(\beta_m^n) g(\beta_m'^n) \xi_m^n,$$

where $\xi_m^n$, $\xi_m'^n$ are given by (6.2), and $\nabla h$ denotes the first derivative of $h$. In fact, it is quite complicated to show (6.16) (especially when $r$ or $l \in (0, 1]$), but it can be proven exactly as in Barndorff-Nielsen, Graversen, Jacod, Podolskij and Shephard [5]. Note also that when $r$ or $l \in (0, 1]$ the terms $\nabla g(\beta_m'^n)$ and $\nabla f(\beta_m^n)$ are still well defined (almost surely), because $\sigma_s \neq 0$ for all $s$. Assumption (V) implies the decomposition

$$\xi_m^n = \xi_m^n(1) + \xi_m^n(2),$$

where $\xi_m^n(1)$, $\xi_m^n(2)$ are defined by

$$\xi_m^n(1) = \frac{n^{1/4}}{n/M - K + 1} \sum_{i=n(m-1)/M}^{n/M-K} \left( \int_{i/n}^{(i+K)/n} (a_u - a_{(m-1)/M}) \, du \right.$$
$$+ \int_{i/n}^{(i+K)/n} \left( \int_{i/n}^{u} a_s' \, ds + \int_{i/n}^{u} (\sigma_{s-}' - \sigma_{(m-1)/M}') \, dW_s \right.$$
$$+ \left. \left. \int_{i/n}^{u} (v_{s-}' - v_{(m-1)/M}') \, dV_s \right) dW_u \right),$$

$$\xi_m^n(2) = \frac{n^{1/4}}{n/M - K + 1} \sum_{i=n(m-1)/M}^{n/M-K} \left( \frac{K}{n} a_{(m-1)/M} + \sigma_{(m-1)/M}' \int_{i/n}^{(i+K)/n} (W_u - W_{i/n}) \, dW_u \right.$$
$$+ \left. v_{(m-1)/M}' \int_{i/n}^{(i+K)/n} (V_u - V_{i/n}) \, dW_u \right),$$

and a similar representation holds for $\xi_m'^n$. Let us now prove that

$$\sum_{m=1}^{M} E[\vartheta_m^n(2)|\mathcal{F}_{(m-1)/M}] \xrightarrow{P} 0. \tag{6.17}$$

A straightforward application of Burkholder's inequality shows that

$$n^{-1/4} \sum_{m=1}^{M} E[\nabla f(\beta_m^n) g(\beta_m'^n) \xi_m^n(1)|\mathcal{F}_{(m-1)/M}] \xrightarrow{P} 0.$$

Next, note that since $f$ is an even function $\nabla f$ is odd. Consequently, $\nabla f(\beta_m^n) g(\beta_m'^n) \xi_m^n(2)$ is an odd functional of $(W, V, B)$. Since $(W, V, B) \stackrel{\mathcal{D}}{=} -(W, V, B)$ we obtain

$$n^{-1/4} \sum_{m=1}^{M} E[\nabla f(\beta_m^n) g(\beta_m'^n) \xi_m^n(2)|\mathcal{F}_{(m-1)/M}] = 0,$$

which implies (6.17). Similarly we can show that

$$\sum_{m=1}^{M} E[\vartheta_m^n(1)|\mathcal{F}_{(m-1)/M}] \xrightarrow{P} 0,$$

which completes the proof of Theorem 4. □

**Proof of Theorem 5.** Theorem 5 can be proven by the same methods as Theorem 4. □

# Acknowledgements

We thank Ole E. Barndorff-Nielsen, Kim Christensen, Holger Dette, Jean Jacod and Neil Shephard for helpful comments and suggestions. The work of the first author was supported by CREATES and funded by the Danish National Research Foundation. The second author is grateful for financial assistance from Deutsche Forschungsgemeinschaft through SFB 475 "Komplexitätsreduktion in multivariaten Datenstrukturen".